\documentclass[10pt,a4paper,twoside]{amsart}

\usepackage[utf8]{inputenc}
\usepackage[english]{babel}
\usepackage[T1]{fontenc}
\usepackage{microtype, fancyhdr, lmodern, xcolor}
\usepackage{listings}

\usepackage[a4paper, hmarginratio=1:1]{geometry}

\usepackage{kerkis}
\usepackage{charter}

\usepackage{amsfonts, amsmath, amsthm, amssymb, mathrsfs, amscd, mathtools}
\usepackage[all]{xy}
\numberwithin{equation}{section}	
\usepackage{faktor}
\allowdisplaybreaks

\usepackage{url, enumitem}
\usepackage[noadjust]{cite}

\usepackage{multirow,bigdelim, caption}
\usepackage{longtable}
\usepackage{booktabs}
\usepackage{graphicx}
\graphicspath{{Figures/}}

\usepackage{hyperref}

\theoremstyle{plain}
\newtheorem{theorem}{Theorem}[section]
\newtheorem{proposition}[theorem]{Proposition}

\theoremstyle{remark}
\newtheorem{remark}[theorem]{Remark}

\theoremstyle{definition}

\newcommand\Zz{\mathbb{Z}}

\newcommand\BB[1]{B_{#1}}	
\newcommand\BC[2]{B_{#1}[#2]}	

\newcommand\sig[1]{\sigma_{\hspace{-0.3ex}#1}^{\null}} 
\newcommand\sigg[2]{\sigma_{\hspace{-0.3ex}#1}^{#2}}	

\newcommand\Sp[2]{\mathrm{Sp}_{#1}(#2)}	

\newcommand\ii{i} 
\newcommand\nn{n} 
\newcommand\nno{n-1} 
\newcommand\mm{m} 

\begin{document}
\title{Powers of half-twists and congruence subgroups of braid groups}

\author[Bellingeri]{Paolo Bellingeri}
\address{Normandie Univ, UNICAEN, CNRS, LMNO, 14000 Caen, France}
\email{paolo.bellingeri@unicaen.fr}

\author[Damiani]{Celeste Damiani}
\address{Fondazione Istituto Italiano di Tecnologia, Genova, Italy}
\email{celeste.damiani@iit.it}

\author[Ocampo]{Oscar Ocampo}
\address{Universidade Federal da Bahia, Departamento de Matem\'atica - IME, CEP:~40170-110 - Salvador, Brazil}
\email{oscaro@ufba.br}

\author[Stylianakis]{Charalampos Stylianakis}
\address{University of the Aegean, Department of mathematics, Karlovasi, 83200, Samos, Greece}
\email{stylianakisy2009@gmail.com}



\date{\today}

\begin{abstract}
In this work, we study the relationship between congruence subgroups $B_n[m]$ and $\mathcal{N}_n(\sigma_1^m)$ the normal closure of $\sigma_1^m$, where $\sigma_1$ is the classical generator of $B_n$. We characterize the conditions under which $\mathcal{N}_n(\sigma_1^m)$ has finite index in $B_n[m]$ and provide explicit generators for these finite quotients. For the cases where the index is infinite, we show that $B_n[m]/\mathcal{N}_n(\sigma_1^m)$ contains a free subgroup. Additionally, we compute the Abelianisation of Coxeter braid subgroups in the finite index cases and construct new finite quotients using commutators of congruence subgroups.

\end{abstract}

\maketitle

\section{Introduction}

Let $n \geq 2$, the \emph{braid group} $\BB\nn$ can be defined as the mapping class group of the $n$-punctured disk, and it admits the following presentation~\cite{Artin:1925}:
\begin{equation}\label{eq:presbn}
\bigg\langle \sig1, \ldots\, , \sig\nno \ \bigg\vert \ 
\begin{matrix}
\sig{i} \sig j = \sig j \sig{i} 
&\text{for} &\vert  i-j\vert > 1\\
\sig{i} \sig j \sig{i} = \sig j \sig{i} \sig j 
&\text{for} &\vert  i-j\vert  = 1
\end{matrix}
\ \bigg\rangle.
\end{equation}

The reducible Burau representation $b\colon \BB\nn \to \mathrm{GL}_n(\mathbb{Z}[t^{\pm 1}])$ is defined by:
\[
b(\sig\ii) = I_{i-1} \oplus \begin{pmatrix}
1-t&t\\	
1&0
\end{pmatrix} \oplus I_{n-i-1},
\]
where $I_k$ is the identity matrix of dimension $k$.

Its specialization at $t=-1$ yields a symplectic representation~\cite[Proposition 2.1]{Gambaudo-Ghys:2005}:
\begin{equation*}
 \rho\colon \BB\nn \rightarrow 
\begin{cases} 
\Sp{\nn-1}{\Zz} \mbox{ for } \nn \mbox{ odd}, \\
(\Sp{\nn}{\Zz})_u \mbox{ for } \nn \mbox{ even}, \\
\end{cases}
\end{equation*}
where $\Sp{\nn}{\Zz}$ denotes the symplectic group  of degree $n$ over $\Zz$ and   $(\Sp{\nn}{\Zz})_u$ is a subgroup of $\Sp{\nn}{\Zz}$ fixing one vector. Via reduction modulo $m$, we obtain:
\begin{equation*}
\label{eqn:rhom}
 \rho_m\colon \BB\nn \rightarrow 
\begin{cases} 
\Sp{\nn-1}{\Zz/m\Zz} \mbox{ for } \nn \mbox{ odd}, \\
(\Sp{\nn}{\Zz/m\Zz})_u \mbox{ for } \nn \mbox{ even}. \\
\end{cases}
\end{equation*}

The kernel of $\rho_m$ ($m \geq 2$) is called the 
\emph{level $m$ congruence subgroup} of $\BB\nn$, denoted $\BC\nn{m}$. For $m=1$, $\BC\nn{1}$ is the \emph{Braid Torelli group}, 
typically denoted $\mathcal{BI}_n$. The above construction is detailed in \cite[Section 2.2]{BDOS:2024}.

Recall that the pure braid group $P_n$ is the kernel of the map from $B_n$ to $S_n$ 
sending each braid to its induced permutation. 
By definition, $\sigg{i}{2}$ belongs to $P_n$ for all $1\leq i\leq n-1$, 
as does the full twist $\Delta_n^2 = (\sig1\sig2 \dots \sig\nno)^n$.
A classical result states that $P_n$ coincides with both the congruence subgroup $B_n[2]$ and $\mathcal{N}_n(\sigma_1^2)$, the normal closure of $\sigma_1^2$ in $B_n$~\cite{Arnold:1968, Brendle-Margalit:2018}. Motivated from the case $m=2$ studied in \cite{Goncalves-Guaschi-Ocampo:2017}, in two previous papers the authors explored the relations between congruence and crystallographic braid groups~\cite{BDOS:2024, BDOS2:2024}, studying the quotient group $\faktor{B_n}{[B_n[m], B_n[m]]}$, where $[B_n[m], B_n[m]]$ is the commutator subgroup of the congruence subgroup $B_n[m]$.

For $m\geq 3$, Coxeter determined the index in $\BB\nn$ of $\mathcal{N}_n(\sigg{1}{m})$, the normal closure of ${\sig1}^m$~\cite{Coxeter:1957}. In the following we denote these normal closures as \emph{Coxeter braid subgroups}. 
It follows from~\cite[Lemma 2.1]{BDOS:2024} that $\mathcal{N}_n(\sigg{1}{m})$ is also a normal subgroup of $\BC\nn\mm$, for any $m\ge 2$. 
Notice that it is natural to study Coxeter-type quotients in braid theory, initially analised for braid groups over the disk by Coxeter \cite{Coxeter:1957} and recently for surface braid groups \cite{DOS}.

In this work we will focus on the relationship between congruence and Coxeter braid subgroups.  
The paper is organised as follows: in Section~\ref{S:Powers} we explicitly compute the quotient  $B_n[m] / \mathcal{N}_n(\sigma_1^m)$ when it is finite. More precisely  in \textbf{Theorem~\ref{main}}
we prove that when $m,n$ are both greater than 2 (the cases $m=2$ or $n=2$ are well known and recalled  in the text) the  quotient  $B_n[m] / \mathcal{N}_n(\sigma_1^m)$ is finite if and only if 
$(m-2)(n-2)<4$ and we determine these finite groups, providing  possible representatives for generators (they are all finite cyclic groups). On the other hand in \textbf{Proposition~\ref{prop:infinite}} we prove that when $n \geq 6$ and $m \notin \{1, 2, 3, 4, 6, 10\}$  the quotient $B_n[m] / \mathcal{N}_n(\sigma_1^m)$ contains a free subgroup.

In Section~\ref{SS:other_quotients} we explore another quotient of $B_n$ involving congruence Coxeter braid subgroups. 
We obtain some finite quotients of $B_n$, in some cases where the Coxeter quotients $B_n / \mathcal{N}_n(\sigma_1^m)$ are infinite. 
More precisely, for positive integers $n$, $m$ and $p$,
we will denote  $\faktor{B_n}{[B_n[m], B_n[m]]}(p)$ the Coxeter-type quotient $\faktor{B_n}{[B_n[m], B_n[m]]\cup \mathcal{N}_n(\sigma_1^p)}$.
The case $m=2$ was studied in~\cite[Theorem~9(a)]{DOS}. 
In \textbf{Proposition 
\ref{prop:crystquot}} we  prove that these groups are finite cyclic groups  when $p=mk+1$ and $k\ge1$ and in \textbf{Proposition 
\ref{prop:crystquot3}} we show that they are finite groups when $m=n=3$ for any integer number $p\ge 3$. 

In Section~\ref{S:abelianizations} we provide (using GAP) the Abelianisation of $\mathcal{N}_n(\sigma_1^m)$  when $(m-2)(n-2)<4$.

\subsection*{Acknowledgments}
P.~B. was partially supported by the ANR project AlMaRe (ANR-19-CE40-0001). C.~D.  is a member of GNSAGA of INdAM.
O.~O. would like to thank Laboratoire de Math\'ematiques Nicolas Oresme (Universit\'e de Caen Normandie) for their hospitality from August 2023 to January 2024, where part of this project was developed and was partially supported by Capes/Programa Capes-Print/ Processo n\'umero 88887.835402/2023-00  and also by National Council for Scientific and Technological Development - CNPq through a \textit{Bolsa de Produtividade} 305422/2022-7. O.~O. also would like to thank Eliane Santos, all HCA staff, Bruno Noronha, Luciano Macedo, Marcio Isabella, Andreia de Oliveira Rocha, Andreia Gracielle Santana, Ednice de Souza Santos, and Vinicius Aiala for their valuable help since July  2024.

\section{Powers of half-twists versus congruence subgroups of braid groups}
\label{S:Powers}

Let us compare the congruence subgroup $\BC{\nn}{\mm}$ and the Coxeter braid subgroup 
$\mathcal{N}_n(\sigg{1}{m})$. For integers $n, m \geq 2$, since $\sigma_1^m \in Ker(\rho_m)=B_n[m]$ (\cite[Lemma 2.1]{BDOS:2024}) and since $Ker(\rho_m)$ is a normal subgroup, the inclusion $\mathcal{N}_n(\sigg{1}{m}) \subseteq B_n[m]$ holds.

Despite this inclusion, ${\mathcal{N}_n(\sigma_1^m)}$ differs significantly from a congruence subgroup. While $\faktor{\BB\nn}{{ B_n[m]}}$ is always finite 
(and partially characterised in \cite{BPS24}), the quotient $\faktor{\BB\nn}{{\mathcal{N}_n(\sigma_1^m)}}$ 
is typically infinite.

Coxeter \cite[Section 10]{Coxeter:1957} characterised when these quotients are finite: 
$\faktor{\BB\nn}{\mathcal{N}_n(\sigma_1^m)}$ is finite if and only if $(m-2)(n-2)<4$ (see also \cite[Chapter~5, Proposition~2.2]{Murasugi}). This inequality yields exactly 
seven cases  where $\mathcal{N}_n(\sigma_1^m)$ has finite index in $\BB\nn$: the families  $(2,m)$ and $(n,2)$ and the five pairs
$(3,3)$, $(4,3)$, $(5,3)$, $(3,4)$, and $(3,5)$. The structure of these quotients is described in Table~\ref{table1_}, where $Q_8$ denotes the quaternion group. These descriptions can be worked out from \cite[Chapter~5]{Murasugi} and can be verified computationally using GAP~\cite{GAP}.

\begin{table}[htb]  
\centering
\begin{tabular}{c || c|c |c}
\hline
Pair $$(n,m)$$ & Quotient group & Order & Structure Description  \\ \hline
 (3,3) & $\faktor{\BB3}{\mathcal{N}_3(\sigma_1^3)}$ & 24 & $SL_2(\Zz/3\Zz)$ \\ 
 (3,4) & $\faktor{\BB3}{\mathcal{N}_3(\sigma_1^4)}$ & 96 & $SL_2(\Zz/3\Zz)\rtimes (\Zz/4\Zz)$ \\ 
 (3,5) & $\faktor{\BB3}{\mathcal{N}_3(\sigma_1^5)}$ & 600 & $\Zz/5\Zz\times SL_2(\Zz/5\Zz)$  \\
 (4,3) & $\faktor{\BB4}{\mathcal{N}_4(\sigma_1^3)}$ & 648 & $((((\Zz/3\Zz)\times (\Zz/3\Zz))\rtimes (\Zz/3\Zz))\rtimes Q_8)\rtimes (\Zz/3\Zz)$\\ 
 (5,3) & $\faktor{\BB5}{\mathcal{N}_5(\sigma_1^3)}$ & 155520 & Extension of $\Zz/6\Zz$ by $PSp_4(\Zz/3\Zz)$ \\ 
\end{tabular}
\caption{Finite quotients of $B_n$ obtained by Coxeter \cite{Coxeter:1957}}
\label{table1_}
\end{table}

Let us examine some known cases. When $n=2$, 
since $B_2 \cong \Zz$ is generated by $\sig1$,
we obtain $\faktor{\BB2}{\mathcal{N}_2(\sigma_1^m)} \cong \Zz / m \Zz$ and that $\mathcal{N}_2(\sigma_1^m) = \BC2{m}$ for any $m\ge 1$.
For $n \geq 3$ and $m=2$, two key facts hold:
$\faktor{\BB\nn}{\mathcal{N}_n(\sigma_1^2)}$ is isomorphic to the symmetric group $S_n$, and
$\mathcal{N}_n(\sigma_1^2)$ coincides with the pure braid group $P_n$, so by Arnol'd's work, we have $\mathcal{N}_n(\sigma_1^2) = \BC{\nn}{2}$ \cite{Arnold:1968, Brendle-Margalit:2018}. 

The equality $\mathcal{N}_3(\sigma_1^3) = \BC3{3}$ holds too, as $\BC3{3}$ is generated by four conjugates of $\sigma_1^3$ \cite[Theorem 4.5]{Stylianakis:2018}.
In what follows, we focus on the Coxeter-type quotient
$\faktor{\BC\nn{m}}{\mathcal{N}_n(\sigma_1^m)}$
for $n,m \geq 3$ and $(m-2)(n-2)<4$, motivated from the Coxeter quotient $\faktor{B_n}{\mathcal{N}_n(\sigma_1^m)}$.

Theorem~\ref{main} is motivated by our prior investigations:

in \cite{BDOS:2024} we raised a question about the quotients of the congruence subgroups by the Coxeter braid subgroups $\BC\nn{m}/\mathcal{N}_n(\sigma_1^m)$.
Here we present a solution. 
Subsection~\ref{SS:other_quotients} further explores the connection between the problem solved in Theorem~\ref{main} and  quotients of braid groups involving commutators of congruence subgroups.

We recall that $\Delta_n^2$ is the Dehn twist (or a full twist) along a curve surrounding all marked points of $D_n$. 
In terms of half-twists it is
$ \Delta_n^2 = (\sig1 \sig2 \ldots \sig\nno)^n$. 
We notice that the full twist $\Delta_n^2$ generates the center of $B_n$ \cite{Chow} and has this notation since it is the square of the Garside element $\Delta_n$, which is another crucial element in braid theory. 
From \cite[Lemma~2.3]{BDOS:2024}, for any $n$, the element $\rho_m(\Delta^2_n)$ is non trivial and of finite order.
Thus, $\Delta^2_n \notin \BC\nn{m}$.

\begin{theorem}
\label{main}
Let $m,n\geq 3$. 
\begin{enumerate}[label=\roman*)]
   \item \label{thm:finite-index}  The Coxeter braid subgroup 
$\mathcal{N}_n(\sigma_1^m)$ has finite index in the congruence subgroup $B_n[m]$ if and only if $(m~-~2)(n~-~2)<~4$.
    \item \label{thm:quotients} In the finite cases $(n,m) \in \{ (3,3), (3,4), (4,3), (3,5), (5,3)\}$, the groups $\faktor{\BC\nn{m}}{\mathcal{N}_n(\sigma_1^m)}$ are as follows. 
    \begin{enumerate}[label=\alph*)]    
        \item For $(n,m)=(3,3)$ and $(4,3)$: $\faktor{\BC\nn{m}}{\mathcal{N}_n(\sigma_1^m)}$ is trivial
        \item For $(n,m)=(3,4)$: $\faktor{\BC3{4}}{\mathcal{N}_3(\sigma_1^4)} \cong \Zz/2\Zz$, 
        generated by $\overline{[\sigma_1^2, \sigma_2^2]}$        \item For $(n,m)=(5,3)$: $\faktor{\BC5{3}}{\mathcal{N}_5(\sigma_1^3)} \cong \Zz/3\Zz$, 
        generated by $\overline{(\sig1 \sig2 \sig3)^4 A \sigma_1^{-2} A^{-1}}$
        where $A=\sigma_4^{-1} \sigma_3^{-2} \sigma_4^{-1} \sigma_2^{-1} \sigma_3 \sigma_2^{-1}$
        \item For $(n,m)=(3,5)$: $\faktor{\BC3{5}}{\mathcal{N}_3(\sigma_1^5)} \cong \Zz/5\Zz$, 
        generated by $\overline{\Delta_3^4}$

    \end{enumerate}
\end{enumerate}
\end{theorem}

\proof 
Let us first prove \ref{thm:finite-index}. 
By \cite[Lemma 2.1]{BDOS:2024} we have that $\mathcal{N}_n(\sigma_1^m) \trianglelefteq \BC\nn{\mm}$. 

Consider the short exact sequence:
\begin{equation}
\label{eq:main}
 1 \to \faktor{\BC\nn{m}}{\mathcal{N}_n(\sigma_1^m)} \to \faktor{\BB\nn}{\mathcal{N}_n(\sigma_1^m)} \to \faktor{\BB\nn}{\BC\nn{m}} \to 1. 
\end{equation}

By Coxeter \cite{Coxeter:1957}, $\faktor{\BB\nn}{\mathcal{N}_n(\sigma_1^m)}$ is finite if and 
only if $(m-2)(n-2)<4$. Since $\faktor{\BB\nn}{\BC\nn{m}}$ is finite, the exact sequence 
implies that $\faktor{\BC\nn{m}}{\mathcal{N}_n(\sigma_1^m)}$ is finite under the same condition. Therefore, $\mathcal{N}_n(\sigma_1^m)$ has finite index in $B_n[m]$ if and only if $(m-2)(n-2)<4$.

\medskip

We now prove \ref{thm:quotients}.
Let $(n,m) \in \{ (3,3), (4,3), (5,3),  (3,5), (3,4)\}$. 

    \begin{enumerate}[label=\alph*)] 

\item As recalled before, the  case  $(3,3)$ is a simple consequence of the fact as $\BC3{3}$ is generated by four conjugates of $\sigma_1^3$ \cite[Theorem 4.5]{Stylianakis:2018}, which implies the other inclusion $\BC3{3}\trianglelefteq \mathcal{N}_3(\sigma_1^3) $. This result can be also easily and directly proved. Consider the short exact sequence \eqref{eq:main} when $n=m=3$:
\[1 \to \faktor{\BC3{3}}{\mathcal{N}_3(\sigma_1^3)} \to \faktor{\BB3}{\mathcal{N}_3(\sigma_1^3)} \to \faktor{\BB3}{\BC\nn{3}} \to 1.\]
By Table \ref{table1_}, $\faktor{\BB3}{\mathcal{N}_3(\sigma_1^3)} \cong SL_2(\Zz/3\Zz)$,  the special linear group of $(2,2)$-matrices with determinant $1$ over $\Zz/3\Zz$, which in turn is isomorphic to $Sp_2(\Zz/3\Zz)$. Also, by definition we have that $\faktor{B_3}{B_3[3]} \cong Sp_2(\Zz/3\Zz)$. Consequently, $\faktor{\BC3{3}}{\mathcal{N}_3(\sigma_1^3)}$ is trivial which implies that $\BC3{3}\cong \mathcal{N}_3(\sigma_1^3)$.

\medskip

 The case $(4,3)$ was already treated by  \cite{Wajnryb:1991}, who provided explicit generators for the kernel of the symplectic representation. However, as for $(3,3)$  we can provide a direct proof. 
 In fact $\faktor{\BB4}{\mathcal{N}_4(\sigma_1^3)}$ is isomorphic to the central extension of the Hessian group by $\Zz/3\Zz$ \cite[Section~4]{Coxeter:1957}. We note that the center of $\faktor{\BB4}{\mathcal{N}_4(\sigma_1^3)}$ is the image of $\Delta_4^2 \in \BB4$ in $\faktor{\BB4}{\mathcal{N}_4(\sigma_1^3)}$.
By \cite[Proposition~4.1]{ArtMich:2006}, the Hessian group is isomorphic to $(\Zz / 3\Zz)^2 \rtimes SL_2(\Zz / 3\Zz)$. On the other hand, combining \cite[Lemmas~2.1, ~2.2, and ~2.3]{BPS24} and specializing to the four strand case we get a short exact sequence
\[ 1 \to \Zz \to (\Sp{4}{\Zz})_u \to \Zz^2 \rtimes \Sp{2}{\Zz} \to 1. \]
Taking the $\mathrm{mod}(3)$ reduction of the short exact sequence above, we derive
\[ 1 \to \Zz / 3\Zz \to (\Sp{4}{\Zz/3\Zz})_u \to (\Zz/3\Zz)^2 \rtimes \Sp{2}{\Zz/3\Zz} \to 1. \]
The abelian group $\Zz / 3\Zz$ from the left hand side of the short exact sequence above is central in $(\Sp{4}{\Zz/3\Zz})_u$ \cite[Lemma~2.4(ii)]{Wajnryb:1991}. Furthermore, $(\Sp{4}{\Zz/3\Zz})_u$ is isomorphic to $\faktor{\BB4}{\BC{4}{3}}$ \cite[Theorem A]{BPS24}.
Since $\mathcal{N}_4(\sigma_1^3) < \BC{4}{3}$, we have a homomorphism $\faktor{\BB4}{\mathcal{N}_4(\sigma_1^3)} \to \faktor{\BB4}{\BC{4}{3}}$, where the center of $\faktor{\BB4}{\mathcal{N}_4(\sigma_1^3)}$ is mapped to the center of $\faktor{\BB4}{\BC{4}{3}}$. Summarizing all the above, we get a commutative diagram.

\begin{equation*}
\xymatrix{
& 1 \ar[r] & \Zz / 3\Zz \ar[r] \ar[d] & \faktor{\BB4}{\mathcal{N}_4(\sigma_1^3)} \ar[r] \ar[d] & (\Zz/3\Zz)^2 \rtimes \Sp{2}{\Zz/3\Zz} \ar[r] \ar[d] & 1 \\
& 1 \ar[r] & \Zz / 3\Zz \ar[r] & \faktor{\BB4}{\BC4{3}} \ar[r] & (\Zz/3\Zz)^2 \rtimes \Sp{2}{\Zz/3\Zz} \ar[r] & 1.
}
\end{equation*}
By the five lemma we get $\faktor{\BB4}{\BC4{3}} \cong \faktor{\BB4}{\mathcal{N}_4(\sigma_1^3)}$ which consequently leads to $\BC4{3} = \mathcal{N}_4(\sigma_1^3)$.

\item We now consider the case $(n,m)=(3,4)$. 
By \cite[Remark~3.7]{BDOS2:2024}, the quotient $\faktor{\BB3}{\BC3{4}}$ is isomorphic to $A_4\rtimes \Zz/4\Zz$, where $A_4$ denotes the alternating group on 4 elements.

By Table~\ref{table1_}, the quotient $\faktor{\BB3}{\mathcal{N}_3(\sigma_1^4)}$ is isomorphic to $\mathrm{SL}_2(\Zz/3\Zz)\rtimes \Zz/4\Zz$. 
We consider once again the short exact sequence given by \eqref{eq:main}:
\[ 1 \to \faktor{\BC3{4}}{\mathcal{N}_3(\sigma_1^4)} \to \faktor{\BB3}{\mathcal{N}_3(\sigma_1^4)} \to \faktor{\BB3}{\BC3{4}} \to 1. \]
The epimorphism $\faktor{\BB3}{\mathcal{N}_3(\sigma_1^4)} \to \faktor{\BB3}{\BC3{4}}$ induces an epimorphism 
\[\mathrm{SL}_2(\Zz/3\Zz)\rtimes \Zz/4  \Zz\to A_4\rtimes \Zz/4\Zz\]
which arises by taking the quotient of $\mathrm{SL}_2(\Zz/3\Zz)$ by its centre.
%
The center of $\mathrm{SL}_2(\Zz/3\Zz)$ has order~2.
Since the quotient of $\mathrm{SL}_2(\Zz/3\Zz)$ by its center 
is isomorphic to $\mathrm{PSL}_2(\Zz/3\Zz)$, and 
$\mathrm{PSL}_2(\Zz/3\Zz) \cong A_4$, 
the induced epimorphism respects the structure of the semidirect product.

By \cite[Proposition~3.3]{BDOS2:2024} or 
\cite[Theorem~1.3]{Nakamura:2021}, $\BC3{4}$ is normally generated by $\sigma_1^4$ and $[\sigma_1^2, \sigma_2^2]$. 
Hence,  that $\faktor{\BC3{4}}{\mathcal{N}_3(\sigma_1^4)}$ is isomorphic to $\Zz/2\Zz$ and it by the image of $[\sigma_1^2, \sigma_2^2]$, which is non-trivial in 
$\faktor{\BC3{4}}{\mathcal{N}_3(\sigma_1^4)}$, and is denoted by$\overline{[\sigma_1^2, \sigma_2^2]}$.

\item For $(n,m)=(5,3)$, let us consider again the short exact sequence~\eqref{eq:main}, that in this case specialises to:
\[ 1 \to \faktor{\BC5{3}}{\mathcal{N}_5(\sigma_1^3)} \to \faktor{\BB5}{\mathcal{N}_5(\sigma_1^3)} \to \faktor{\BB5}{\BC5{3}} \to 1. \]

By definition $\faktor{\BB5}{\BC5{3}}$ lies inside $\Sp4{\Zz/3\Zz}$. By ~\cite[Theorem 1 (1)]{ACampo:1979} the latter inclusion is an isomorphism. Also, $\Sp4{\Zz/3\Zz}$ has order 51840. From Table~\ref{table1_}, the quotient $\faktor{\BB5}{\mathcal{N}_5(\sigma_1^3)}$ has order 155520. Hence, $\faktor{\BC5{3}}{\mathcal{N}_5(\sigma_1^3)}$ 
is isomorphic to~$ \Zz/3\Zz$, 
the cyclic group of 
order 3.

 By \cite[Theorem 4.4]{Stylianakis:2018} 
 we have that $\BC5{3}$ is normally  generated by $\mathcal{N}_5(\sigma_1^3)$ together with the element $(\sig1 \sig2 \sig3)^4 A \sigma_1^{-2} A^{-1}$ where 
$A=\sigma_4^{-1} \sigma_3^{-2} \sigma_4^{-1} \sigma_2^{-1} \sigma_3 \sigma_2^{-1}$. Therefore the image of $A$ generates the quotient $\faktor{\BC5{3}}{\mathcal{N}_5(\sigma_1^3)}$.

\medskip

\item Finally, let us consider $(n,m)=(3,5)$, consider again the short exact sequence~\eqref{eq:main}:
\[ 1 \to \faktor{\BC3{5}}{\mathcal{N}_3(\sigma_1^5)} \to \faktor{\BB3}{\mathcal{N}_3(\sigma_1^5)} \to \faktor{\BB3}{\BC3{5}} \to 1. \]
We claim that 
$\faktor{\BC3{5}}{\mathcal{N}_3(\sigma_1^5)} \cong \Zz/5\Zz$. To compute the generator of $\faktor{\BC3{5}}{\mathcal{N}_3(\sigma_1^5)}$, we will prove that $\faktor{\BB3}{\BC3{5}} \cong \faktor{\BB3}{\mathcal{N}_3(\sigma_1^5, \Delta_3^4)}$, where $\Delta_3^2 = (\sig1 \sig2)^3$. Therefore, $\faktor{\BC3{5}}{\mathcal{N}_3(\sigma_1^5)}$ is normally generated by the lift of $\Delta_3^2$ in $\faktor{\BB3}{\mathcal{N}_3(\sigma_1^5)}$. Since, $\faktor{\BC3{5}}{\mathcal{N}_3(\sigma_1^5)}$ is a cyclic group, then the latter group is generated by a single element, the image of $\Delta_3^2$ in $\faktor{\BB3}{\mathcal{N}_3(\sigma_1^5)}$.

By definition $\faktor{\BB3}{\BC3{5}}$ lies inside $\Sp2{\Zz/5\Zz}$ and by \cite[Theorem 1 (1)]{ACampo:1979} the latter inclusion is an isomorphism.

By Wajnryb \cite[Theorem 1]{Wajnryb:1991}, 
 we have the presentation:
\[ \Sp{2}{\Zz/5} = \langle \sig1, \sig2 \mid \sig1 \sig2 \sig1 = \sig2 \sig1 \sig2 , \sigg{1}{5} = 1, (\sig1 \sig2)^6 = 1, \sigg{1}{3} \sigg{2}{4} \sigg{1}{2} = \sigg{2}{-2} \sig1 \sigg{2}{2} \rangle. \]

The last relation in this presentation is redundant, as we show below. In the following computation, we use the underlines in order to indicate the manipulated parts.

\begin{align*}
\sigg{1}{-2} \underline{\sigg{2}{-4}} \sigg{1}{-3} \sigg{2}{-2} \sig1 \sigg{2}{2} & = \sigg{1}{-2} \sig2 \sigg{1}{-3} \sigg{2}{-3} \underline{\sig2 \sig1 \sig2} \sig2 \\
& = \sigg{1}{-2} \sig2 \sigg{1}{-3} \underline{\sigg{2}{-3} \sig1 \sig2 \sig1 \sig2}\\
& = \sigg{1}{-2} \sig2 \sigg{1}{-3} \sigg{2}{-4} (\sig2 \sig1 \sig2 \sig1 \sig2)\\
& = \sigg{1}{-2} \sig2 \sigg{1}{-3} \sigg{2}{-4} \sigg{1}{-1} (\sig1 \sig2)^3\\
& = \sigg{1}{-2} \sig2 \sigg{1}{-3} \underline{\sigg{2}{-4}} \sigg{1}{-1} (\sig1 \sig2)^3\\
& = \sigg{1}{-2} \sig2 \sigg{1}{-3} \sig2 \sigg{1}{-1} (\sig1 \sig2)^3 \\
& = \sigg{1}{-2} \sig2 \underline{\sigg{1}{-4}} \sig1 \sig2 \sig1 \sigg{1}{-2} (\sig1 \sig2)^3 \\
& = \sigg{1}{-2} \sig2 \sig1 \underline{\sig1 \sig2 \sig1} \sigg{1}{-2} (\sig1 \sig2)^3 \\
& = \sigg{1}{-2} \sig2 \sig1 \sig2 \sig1 \sig2 \sigg{1}{-2} (\sig1 \sig2)^3 \\
& = \sigg{1}{-3} (\sig1 \sig2 \sig1 \sig2 \sig1 \sig2) \sigg{1}{-2} (\sig1 \sig2)^3 \\
& = \sigg{1}{-3} (\sig1 \sig2)^3 \sigg{1}{-2} (\sig1 \sig2)^3 \\
& = \sigg{1}{-3} \sigg{1}{-2} (\sig1 \sig2)^3 (\sig1 \sig2)^3\\
& = \sigg{1}{-5} (\sig1 \sig2)^6\\
& = 1.
\end{align*}

Therefore
\[ \Sp2{\Zz/5\Zz} \cong \faktor{\BB3}{\mathcal{N}_3(\Delta_3^4, \sigma_1^5)}. \]

Thus, $\faktor{\BC3{5}}{\mathcal{N}_3(\sigma_1^5)}$ is generated by $\overline{\Delta_3^4}$.

\end{enumerate}
\endproof

\begin{remark}
 Theorem \ref{main} has been independently and recently claimed in \cite{Banerjee-Huxford} (Remark~2.1). However, here we give a detailed proof giving also an explicit generator for the different finite quotients.
\end{remark}

When $(m-2)(n-2)>4$ very little is known on the structure of $\faktor{\BC\nn{m}}{\mathcal{N}_n(\sigma_1^m)}$.
Here we can show that under certain conditions on $n$ and~$m$, it contains a copy of a free group. 

\begin{proposition}\label{prop:infinite}
Let $n\geq 6$ and let $m$ be a positive integer such that $m \notin \{1,2,3,4,6,10 \}$. Then, $\faktor{\BC\nn{m}}{\mathcal{N}_n(\sigma_1^m)}$ contains a copy of a free group.
\end{proposition}

\proof
Denote by $M_{0,n}$ the group defined by the following presentation.

\begin{equation*}\label{eq:presMn}
\left\langle  \sig1, \ldots\, , \sig\nno \ \bigg\vert \ 
\begin{matrix}
\sig{i} \sig j = \sig j \sig{i} 
&\text{for} &\vert  i-j\vert > 1\\
\sig{i} \sig j \sig{i} = \sig j \sig{i} \sig j 
&\text{for} &\vert  i-j\vert  = 1\\
(\sig1 \sig2 ... \sig\nno)^n = 1 \\
\sig1 \sig2 \cdots \sigma_{n-1}^2 \cdots \sig1 = 1
\end{matrix}
\ \right\rangle.
\end{equation*}

From the presentation we can see that $M_{0,n}$ is a quotient of $\BB\nn$. Remark that the group $M_{0,n}$ is the mapping class group of a sphere with $n$ marked points (see also~\cite{BirmanHilden:1969} and \cite{Farb-Margalit:book}).
For $f\in M_{0,n}$, we denote by $\overline{\mathcal{N}_n}(f)$ the normal closure of $f$ in $M_{0,n}$. This gives us a homomorphism
\[ \faktor{\BB\nn}{\mathcal{N}_n(\sigma_1^m)} \to \faktor{M_{0,n}}{\overline{\mathcal{N}_n}(\sigma_1^m)}. \]

By \cite[Theorem B]{Stylianakis:2019}, $\faktor{M_{0,n}}{\overline{\mathcal{N}_n}(\sigma_1^m)}$ contains a free group, when $m \notin \{1,2,3,4,6,10 \}$ . Consequently, $\faktor{\BB\nn}{\mathcal{N}_n(\sigma_1^m)}$ contains a free subgroup. 
The result then follows from the short exact sequence~\ref{eq:main}, which shows that $\faktor{\BC\nn{m}}{\mathcal{N}_n(\sigma_1^m)}$ contains a free subgroup.
\endproof

\begin{remark}
    The argument in \cite[Theorem B]{Stylianakis:2019} relies on the Jones representation of
the mapping class group of a punctured sphere and it does not give a method to explicitly calculate the rank of the free group.
\end{remark}

\section{Other related finite quotients}
\label{SS:other_quotients}
In two previous papers the authors explored the relations between congruence and 
crystallographic braid groups~\cite{BDOS:2024, BDOS2:2024}.

A natural question arising from this study concerns the finiteness of further quotients of the group~$\faktor{\BC\nn{m}}{\mathcal{N}_n(\sigma_1^m)}$.
More precisely, for positive integers $n$, $m$ and $p$,
let us denote by $\faktor{B_n}{[B_n[m], B_n[m]]}(p)$ the Coxeter-type quotient of $\faktor{B_n}{[B_n[m], B_n[m]]}$ by adding the additional relations $\sigma_i^p~=~1$, for all $1\leq i\leq n-1$. Our question is: for which values of $n$, $m$, and $p$ is the group $\faktor{B_n}{[B_n[m], B_n[m]]}(p)$ finite?

Some partial answers to this question are known: recalling that $B_n[2]=P_n$, 
 in~\cite[Theorem~9(a)]{DOS}
 it is proved that the Coxeter-type quotient $\faktor{B_n}{[P_n, P_n]}(p)$ is finite for all $n\geq 3$ and all $p$ and its order has been explicitly computed.

Here we can provide a new condition of finiteness when $m>2$ based on a previous result of \cite{BDOS:2024}.

\begin{proposition}
\label{prop:crystquot}
    Let $n\geq 3$, $m\geq 2$ and $k\geq 1$. The group $\faktor{B_n}{[B_n[m], B_n[m]]}(mk+1)$ is finite and isomorphic to the cyclic group $\Zz_{mk+1}$. 
\end{proposition}

\proof
Let $n\geq 3$, $m\geq 2$ and $k\geq 1$. 
By \cite[Lemma~2.1]{BDOS:2024}, we have that
$\mathcal{N}_n(\sigma_1^m)$ is normal in $\BC\nn{\mm}$.

From the definition of the quotient group 
$\faktor{B_n}{[B_n[m], B_n[m]]}(mk+1)$, 
each generator satisfies the relation $\sigma_i^{mk+1}=1$, 
or equivalently $\sigma_i=(\sigma_i^m)^{-k}$. This implies the commutation relation $[\sigma_i, \sigma_{i+1}]=1$, for all $1\leq i\leq n-1$.

The braid relation $\sigma_i\sigma_{i+1}\sigma_i=\sigma_{i+1}\sigma_i\sigma_{i+1}$
together with commutation then implies that $\sigma_i=\sigma_{i+1}$  in the quotient group $\faktor{B_n}{[B_n[m], B_n[m]]}(mk+1)$,
for all $1\leq i\leq n-1$.  

We can then conclude that $\faktor{B_n}{[B_n[m], B_n[m]]}(mk+1)\cong \Zz_{mk+1}$.
\endproof

\begin{remark}
   In the case $m=2$, the last proposition recovers exactly \cite[Theorem~9(a)]{DOS}.
\end{remark}

Now, we consider the special case of the congruence subgroup $B_3[3]$. In \cite[Theorem~3.13]{BDOS2:2024} was proved that  the group $\faktor{B_3}{[B_3[3], B_3[3]]}$ is torsion free, in the following proposition we prove that the corresponding Coxeter-type quotient is finite.  

\begin{proposition}
\label{prop:crystquot3}
    For any integer number $q\geq 3$ the group $\faktor{B_3}{[B_3[3], B_3[3]]}(q)$ is finite. 
\end{proposition}

\proof

Let $q\geq 3$ be an integer number. If $q$ is divided by 3 then from the division algorithm we have that there exist unique integers $k$ and $r$ with $0\leq r \leq 2$ such that $q=3k+r$. 
In the first part of this proof, we examine the group $\faktor{\BC3{3}}{[\BC3{3},\BC3{3}]}(q)$ analising separately the cases $r=0,1,2$. 

If $r=1$ and so $q=3k+1$ the result follows from Proposition~\ref{prop:crystquot}. 

Recall from \cite[Section~3.3]{BDOS2:2024} that the group $\faktor{B_3[3]}{[B_3[3], B_3[3]]}$ is free Abelian of rank~4, and a basis is given by $\{ \sigma_1^3,\, \sigma_2^3,\, \sigma_1\sigma_2^3\sigma_1^{-1},\, \sigma_1^{-1}\sigma_2^3\sigma_1 \}$. 
Now, suppose that $r=0$ and so $q=3k$. Hence $(\sigma_i^3)^{k}=1$, for $i=1,2$.   
If $r=2$ then $q=3k+2$ and in this case we have that $(\sigma_i^3)^{3k+2}=1$, for $i=1,2$.  Then, in any case, the group $\faktor{\BC3{3}}{[\BC3{3},\BC3{3}]}(q)$ is an abelian finite group obtained by taking the quotient of $\Zz^4$ adding the only relations that the generators have finite order.

Finally, consider the short exact sequence that comes from a sequence that was used in \cite[Section~3.3]{BDOS2:2024} to prove that the crystallographic group $\faktor{\BB3}{[\BC3{3},\BC3{3}]}$ is torsion free
\[
\begin{CD}
1 @>>> \faktor{\BC3{3}}{[\BC3{3},\BC3{3}]}(q) @>>> \faktor{\BB3}{[\BC3{3},\BC3{3}]}(q) @>>> \Sp{2}{\Zz/3} @>>> 1.
\end{CD}
\]
Therefore, the middle group $\faktor{\BB3}{[\BC3{3},\BC3{3}]}(q)$ is an extension of two finite groups, proving this result.
\endproof

 \section{Abelianisations}
 \label{S:abelianizations}
 
Our study of quotients $\faktor{\BC\nn{m}}{\mathcal{N}_n(\sigma_1^m)}$ naturally leads to examining the structure of $\mathcal{N}_n(\sigma_1^m)$ itself. In this section we compute its Abelianisation for the finite-index cases.

\begin{theorem}
If $(n,m) \in \{ (3,3), (4,3), (5,3), (3,4), (3,5) \}$, then, $\mathrm{H}_1(\mathcal{N}_n(\sigma^m);\Zz)$ is torsion free of finite rank as follows. 

\[ 
\mathrm{rank}(\mathrm{H}_1(\mathcal{N}_n(\sigma_i^m);\Zz)) = \left\{
\begin{array}{ll}
      4, & \mathrm{if} \quad (n,m) = (3,3) \\
      6, & \mathrm{if} \quad (n,m) = (3,4) \\
      12, & \mathrm{if} \quad (n,m) = (3,5) \\
      12, & \mathrm{if} \quad (n,m) = (4,3) \\
      40, & \mathrm{if} \quad (n,m) = (5,3) \\
\end{array} 
\right. 
\]
\label{rankcoxeter}
\end{theorem}
\proof
If $(n,m)$ is as in Theorem \ref{rankcoxeter}, then $\mathcal{N}_n(\sigma_i^m)$ is a finite index subgroup of $\BB\nn$ (see table \ref{table1_}). Then one can apply the Coxeter-Todd algorithm to get the Schreier transversals of $\mathcal{N}_n(\sigma_i^m)$ in $\BB\nn$ and apply the Reidemeister-Schreier method to compute a presentation of $\mathcal{N}_n(\sigma_i^m)$ \cite[Chapters 8 and 9]{Johnson}. Furthermore, having a presentation of a group, one can compute its abelianization (and consequently its rank) by quotient it by all commutators and applying Tietze transformations \cite[Section 4.4]{Johnson}.

Since the orders of $\faktor{\BB\nn}{\mathcal{N}_n(\sigma_i^m)}$ are large enough for a calculation by hand, so we compute these abelianizations using GAP. However, it is possible to compute the abelianization of $\mathcal{N}_3(\sigma_i^3)$ by a different method in the next Section.

To explain the steps of the algorithm below, we first define a free group whose generators are in 1-1 correspondence with the generators of $\BB\nn$. Then we define presentations of $\BB\nn$ (represented by $g$) and $\mathcal{N}_n(\sigma^m)$ (represented by $g1$). In the next step we define a homomorphism, denoted by $map$, between $g$ and $g1$ (which is actually an epimorphism) and the kernel of $map$. The function \verb|GeneratorsOfGroup| gives the generators of the kernel $k$ and then we use the \verb|AbelianInvariants| function to determine the Abelianisation structure:

\begin{verbatim}
f:=FreeGroup("a","b", "c", ...);
g:=f/ParseRelators(f,"a*b*a*b^-1*a^-1*b^-1, ...");
g1:=f/ParseRelators(f,"a*b*a*b^-1*a^-1*b^-1, ..., a^m");
map:=GroupHomomorphismByImages(g,g1,
    GeneratorsOfGroup(g),GeneratorsOfGroup(g1));
k := Kernel(map);
GeneratorsOfGroup(k);
AbelianInvariants(k);
\end{verbatim}

\endproof

To examine the structure more closely, we analyse the generators and braid action for the case $(n,m)=(3,4)$:

\begin{table}[htb] 
\centering
\begin{tabular}{c|c}
\hline
Generators of $\mathcal{N}_3(\sigma^4)$ & mapping $x_i$  \\ \hline

$\sigma_1^{-4} $ & $x_1$ \\ 
$\sigma_2^{-4} $ & $x_2$ \\ 
$\sigma_1 \sigma_2^{-4} \sigma_1^{-1} $ & $x_3$ \\ 
$\sigma_1^{-1} \sigma_2^{-4} \sigma_1 $ & $x_4$ \\ 
$\sigma_1^2 \sigma_2^{-4} \sigma_1^{-2} $ & $x_5$ \\ 
$\sigma_2 \sigma_1^{-1} \sigma_2^3 \sigma_1^{-1} \sigma_2 \sigma_1 $ & $x_6$

\end{tabular}
\caption{For convenience, in the case $(n,m)=(3,4)$, we rename the generators of $\mathcal{N}_3(\sigma^4)$ to $x_i$.}
\label{table2}
\end{table}

\begin{table}[htb] 
\centering
\begin{tabular}{c|c|c}
\hline
Generators of $\BB3$ & Generators of $\mathcal{N}_3(\sigma^4)$ & Action of $\BB3$ on $\mathcal{N}_3(\sigma^4)$  \\ \hline
$\sigma_1$ & $x_1$ & $\sigma_1 x_1 \sigma_1^{-1}$ = $x_1$ \\ 
$\sigma_1$ & $x_2$ & $\sigma_1 x_2 \sigma_1^{-1}$ = $x_3$ \\ 
$\sigma_1$ & $x_3$ & $\sigma_1 x_3 \sigma_1^{-1}$ = $x_5$ \\ 
$\sigma_1$ & $x_4$ & $\sigma_1 x_4 \sigma_1^{-1}$ = $x_2$ \\ 
$\sigma_1$ & $x_5$ & $\sigma_1 x_5 \sigma_1^{-1}$ = $x_1^{-1} x_4 x_1$ \\ 
$\sigma_1$ & $x_6$ & $\sigma_1 x_6 \sigma_1^{-1}$ = $x_2 x_6 x_2^{-1}$ \\
												
$\sigma_2$ & $x_1$ & $\sigma_2 x_1 \sigma_2^{-1}$ = $x_4$ \\
$\sigma_2$ & $x_2$ & $\sigma_2 x_2 \sigma_2^{-1}$ = $x_2$ \\
$\sigma_2$ & $x_3$ & $\sigma_2 x_3 \sigma_2^{-1}$ = $x_1$ \\
$\sigma_2$ & $x_4$ & $\sigma_2 x_4 \sigma_2^{-1}$ = $x_6^{-1}$ \\
$\sigma_2$ & $x_5$ & $\sigma_2 x_5 \sigma_2^{-1}$ = $x_1 x_5 x_1^{-1}$ \\
$\sigma_2$ & $x_6$ & $\sigma_2 x_6 \sigma_2^{-1}$ = $x_2^{-1} x_3^{-1} x_2$ 
\end{tabular}
\caption{For the case $(n,m)=(3,4)$ we get the action of $\BB3$ on $\mathcal{N}_3(\sigma^4)$. }
\label{table3}
\end{table}

Such actions can be computed using GAP:

\begin{verbatim}
hom:=EpimorphismFromFreeGroup(k:names:=["x1","x2","x3", ...]);
PreImagesRepresentative(hom, g.i*k.j*g.i^-1);
\end{verbatim}
where \verb|hom| maps the kernel generators to the names $x_i$ and \verb|PreImagesRepresentative| factorizes conjugations in terms of these generators.

The pattern generalizes to larger cases: $\mathcal{N}_3(\sigma^5)$ and $\mathcal{N}_4(\sigma^3)$ each have 12 generators, while $\mathcal{N}_5(\sigma^3)$ has 40 generators. The full generators and actions for these cases are not included but they can be recovered using GAP~\cite{GAP}.

\endproof

\begin{remark}
In the case $(n,m)=(3,3)$ we have that $\mathcal{N}_3(\sigma^3) = \BC3{3}$ \cite{Stylianakis:2018}. A special case of a result of the Authors states that $\BC3{3} \cong F_3 \times \mathbb{Z}$, where $F_3$ is a free group of rank $3$ \cite[Theorem 3.9]{BDOS2:2024}. As a consequence, we get that $H_1(\mathcal{N}_3(\sigma^3);\mathbb{Z})=H_1(\BC3{3};\mathbb{Z}) = \mathbb{Z}^4$.
\end{remark}

\bibliography{braid_groups.bib}{}
\bibliographystyle{alpha}
\end{document}